\documentclass[twoside,11pt]{article}

  \usepackage{graphicx}

  \usepackage[latin1]{inputenc} %pour ne pas avoir besoin de coder les accents

  \usepackage{amsmath}
  \usepackage{amssymb}
  \usepackage{amsthm}
  \usepackage{a4}
  \usepackage{mathrsfs}
  \usepackage{bbm}
  \usepackage{color}
  \usepackage{graphicx}
  \usepackage[normal,small,sc,md,sf]{caption2}

\font\MacDo=macdo12
\font\MacDoSmall=macdo7

  \newcommand{\macdo}{\text{\MacDo\char 77}}
\newcommand\ind[1]{\mathbf{1}_{#1}} 
\newcommand{\nb}{\overline{n}}
\newcommand{\kb}{\overline{k}}
\newcommand{\ta}{\text{\MacDo\char 65}}%{{\bigtriangledown}} % Le symbole pour les tableaux triangulaires
\newcommand{\tal}{\text{\MacDo\char 66}}%{{\bigtriangledown^{(L)}}} % Le symbole pour le sous-tableau de gauche
\newcommand{\tar}{\text{\MacDo\char 67}}%{{\bigtriangledown^{(R)}}} % Le symbole pour le sous-tableau de droite

\newcommand{\sta}{\text{\MacDoSmall\char 65}}%{{\bigtriangledown}} % Le symbole pour les tableaux triangulaires
\newcommand{\stal}{\text{\MacDoSmall\char 67}}%{{\bigtriangledown^{(L)}}} % Le symbole pour le sous-tableau de gauche
\renewcommand{\star}{\text{\MacDoSmall\char 66}}%{{\bigtriangledown^{(R)}}} % Le symbole pour le sous-tableau de droite

\newcommand{\sign}{\mathrm{sign}}

\newcommand{\h}{f}

\newcommand{\converge}[2]{\mathop{\longrightarrow}_{\scriptscriptstyle{#1}}^{\scriptscriptstyle{#2}}}
  \renewcommand{\P}{\mathscr{P}}
  \newtheorem{theorem}{Theorem}[section]
  \newtheorem{lemma}{Lemma}[section]
  
  \newtheorem{proposition}{Proposition}[section]
\begin{document}

    \setcaptionmargin{1cm}
    \setlength{\captionindent}{0.5cm}
\title{Self-Similar Corrections to the Ergodic Theorem for the Pascal-Adic Transformation}
\author{\'E. Janvresse, T. de la Rue, Y. Velenik}

\maketitle

\begin{abstract}
 Let $T$ be the Pascal-adic transformation. For any measurable function $g$,
 we consider the corrections to the ergodic theorem
 $$ \sum_{k=0}^{j-1}g(T^kx)-\dfrac{j}{\ell}\sum_{k=0}^{\ell-1}g(T^kx). $$
 When seen as graphs of functions defined on $\{0,\ldots,\ell-1\}$,
 we show for a suitable class of functions $g$ that these quantities, once properly renormalized, converge to
(part of) the graph of a self-affine function. The latter only depends on the ergodic component of $x$, 
and is a deformation of the so-called Blancmange function. We also briefly describe the links with a series of works on Conway recursive \$10,000  sequence.
\end{abstract}



\section{Introduction}

\subsection{The Pascal-adic transformation and its basic blocks}

The notion of \emph{adic transformation} has been introduced by Vershik (see e.g. \cite{versh5}, \cite{versh6}),
as a model in which the transformation acts on infinite paths in some graphs, called 
\emph{Bratteli diagrams}. As shown by Vershik, every ergodic automorphism of the Lebesgue space
is isomorphic to some adic transformation, with a Bratteli diagram which may be quite
complicated. Vershik also proposed to study the ergodic properties of an adic transformation
in a given simple graph, such as the Pascal graph which gives rise to the so-called 
\emph{Pascal adic transformation}. We recall the construction of the latter by the cutting-and-stacking
method in the appendix.

\begin{figure}[ht]
   \centering
   \includegraphics[scale=0.62]{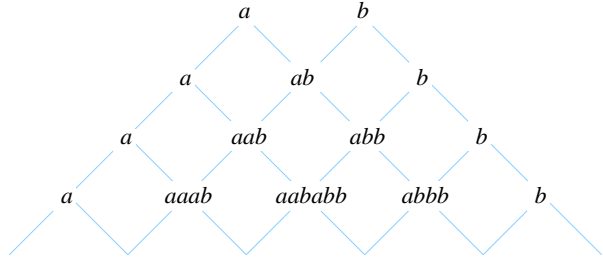}
   \caption{The beginning of the words triangle.}
   \label{fig:words}
\end{figure}

When studying the Pascal-adic transformation, one is naturally
led to consider the family of words $B_{n,k}$ ($n\ge 1,\ 0\le k\le n$)
on the alphabet $\{a,b\}$, inductively defined by (see Figure~\ref{fig:words})
$$ B_{n,0}\ :=\ a,\quad B_{n,n}\ :=\ b,\quad(n\ge 1)$$
and for $0<k<n$
$$ B_{n,k}\ :=\ B_{n-1,k-1}B_{n-1,k}. $$
It follows easily from this definition that the length of the block $B_{n,k}$ is given by the binomial coefficient $\binom{n}{k}$. 

\begin{figure}[t!]
   \centering
   \includegraphics[height=1cm]{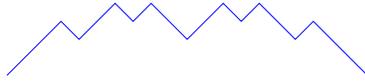}
   \caption{The graph $F_{6,3}$ associated to the word 
$B_{6,3}=aaabaababbaababbabbb$.}
   \label{fig:f63}
\end{figure}

In order to describe the large-scale structure of the basic blocks $B_{n,k}$, we associate to each of them the graph of a real-valued function $F_{n,k}$. Let us denote by $B_{n,k}(\ell)$ the $\ell$th letter of $B_{n,k}$.
For each $n\ge2$ and $0\le k\le n$, we consider the function $F_{n,k}:[0,\binom{n}{k}]\to\mathbb{R}$ defined from the basic
block $B_{n,k}$ as follows (see Figure~\ref{fig:f63}):
\begin{itemize}
\item $F_{n,k}(0)=0$;
\item if $1\le \ell\le \binom{n}{k}$ is an integer, $F_{n,k}(\ell) =
\begin{cases}
 F_{n,k}(\ell-1)  +1 & \text{if $B_{n,k}(\ell) = a$,} \\
 F_{n,k}(\ell-1)  -1 & \text{if $B_{n,k}(\ell) = b$;}
\end{cases}$
\item $F_{n,k}$ is linearly interpolated between two consecutive integers.
\end{itemize}

\begin{figure}[t]
   \centering
   \includegraphics[width=6cm, height=4.5cm]{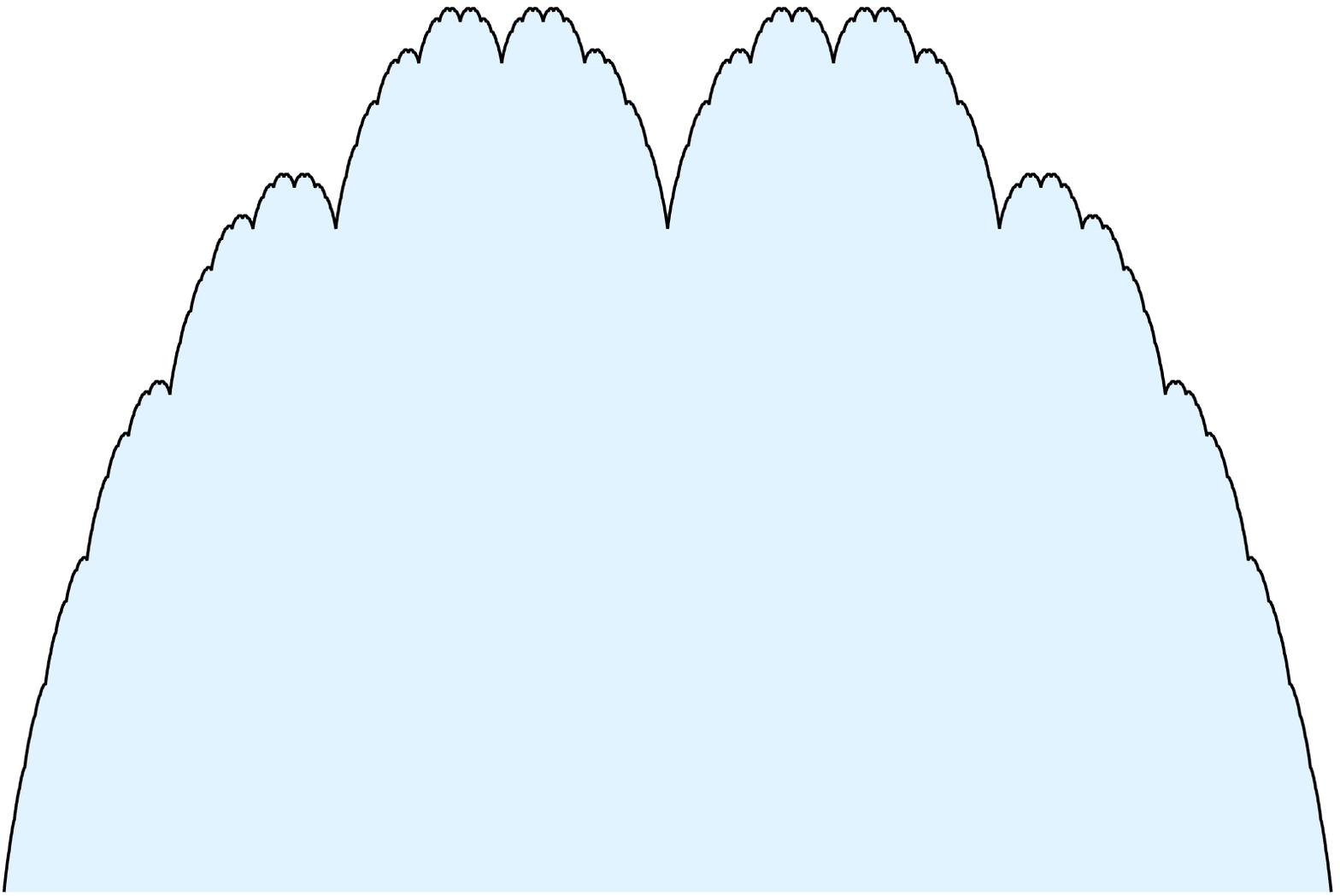}\quad 
   \includegraphics[width=6cm, height=4.5cm]{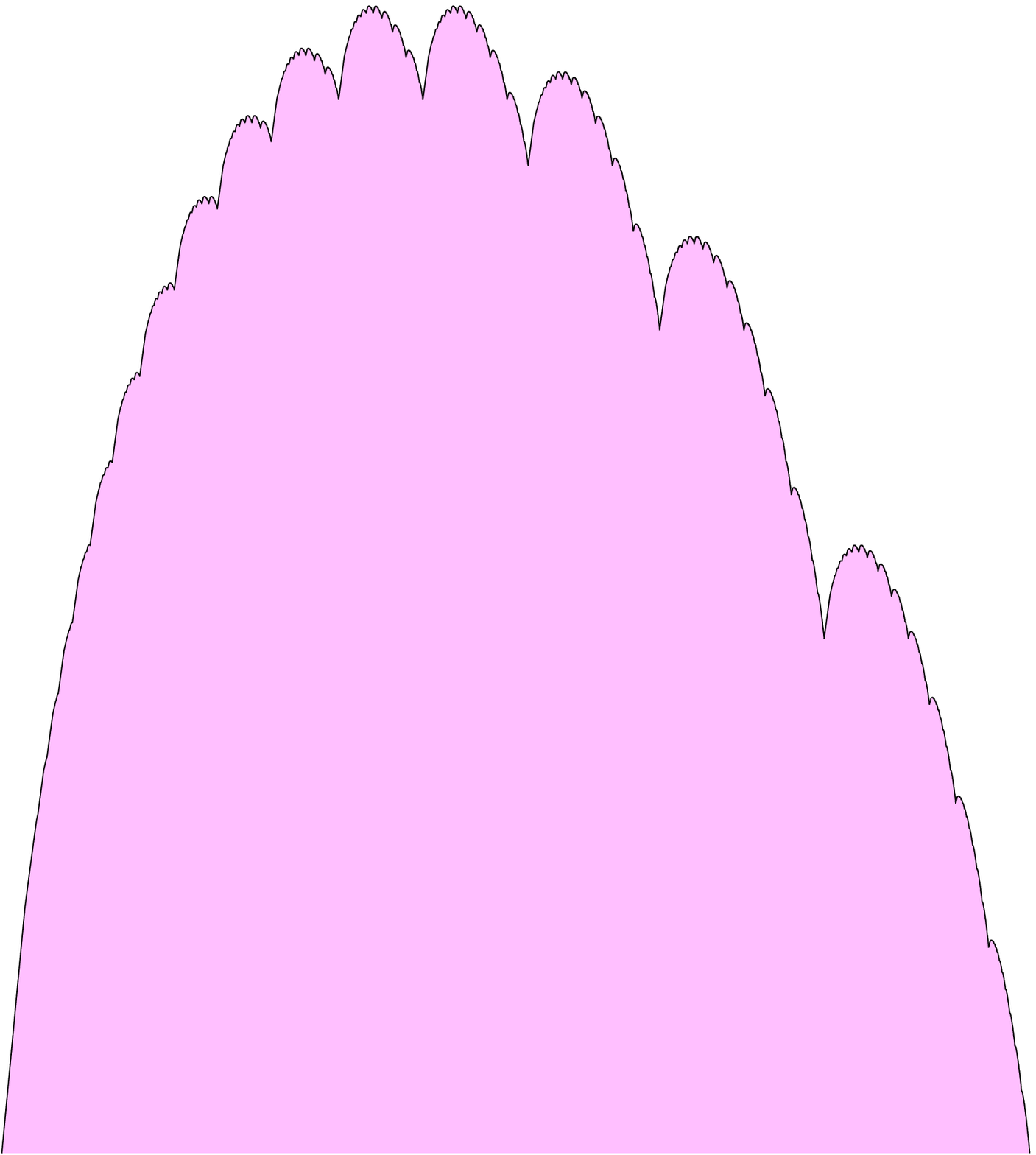}
   \caption{Limiting observed shape for $B_{n,pn}$ with $p=0.5$ (left) and $p=0.8$ (right).}
   \label{fig:macdo0.5}
\end{figure}

As we will see in Section~\ref{ErgodicInterpretation}, the ergodic theorem implies that
the graph of the function
$$
t\ \longmapsto\ 
\dfrac{1}{\tbinom{n}{k}} F_{n,k}\left(t\tbinom{n}{k}\right)
$$
converges to a straight line as $n\to\infty$ and $k/n\to p$. In order to extract the nontrivial structure of 
this graph, we have to remove this dominant contribution and look at the correction (see Section~\ref{limit_standard}). Once this is done, 
it appears that the resulting graph  converges to the graph of a self-affine function depending only on $p=\lim k/n$, described in the following section.
Examples of  such limiting graphs are shown in
Figure~\ref{fig:macdo0.5}. 

\subsection{A one-parameter family of self-affine maps}

\begin{figure}[t]
   \centering
   \includegraphics[scale=.5]{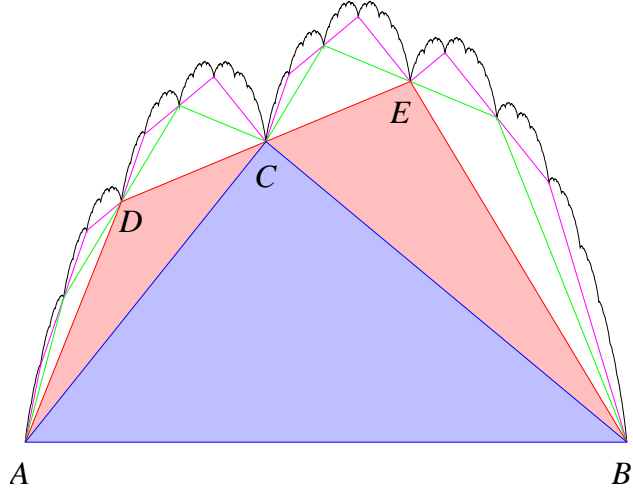}
   \caption{The first four steps in the construction of $\macdo_p$, for $p=0.4$. 
   In the first step, we transform the original interval $AB$ into the polygonal line $ACB$, where $C=\alpha_p^L(B)=\alpha_p^R(A)$. In the second step, we similarly map the segment $AC$ onto the segments $AD$ and $DC$, and the segment $CB$ onto $CE$ and $EB$, by applying the same affine transformations. The procedure is then iterated yielding an increasing sequence of piecewise-linear functions, converging to the self-affine function $\macdo_p$.
   }
   \label{fig_affinite}
\end{figure}

For any $0<p<1$, we consider the two affinities $\alpha_p^L$ and $\alpha_p^R$
defined by
$$ \alpha_p^L(x,y)\ :=\ (px,py+x), $$
and
$$ \alpha_p^R(x,y)\ :=\ \Bigl((1-p)x +p,(1-p)y-x+1\Bigr). $$
These maps are strict contractions of $[0,1]\times \mathbb{R}$,
thus there exists a unique compact set $E_p$ such that
$$ E_p\ =\ \alpha_p^L(E_p)\cup\alpha_p^R(E_p). $$
As shown in \cite{Bedford89}, $E_p$ is the graph of a continuous
self-affine map $\macdo_p\,:\ [0,1]\to \mathbb{R}$, whose construction is
illustrated in Figure~\ref{fig_affinite} (see also \cite[Chapter 11]{Falconer90}).

Note that $\macdo_{0.5}$ is known as the ``Blancmange function'', or ``Takagi fractal curve'', and was introduced in~\cite{Takagi1903}.

\subsection{Conway recursive \$10,000 sequence}

In a lecture at AT\& T Bell Labs in 1988, Conway introduced the following recursive sequence,
$$
C(n) = C(C(n-1)) + C(n-C(n-1)),
$$
with initial conditions $C(1)= C(2)=1$. The latter was then studied and generalized in a large number of papers, see e.g.~\cite{Mallows1991, KuboVakil1996}. In Appendix~\ref{app_conway}, we briefly describe some links between this topic and the content of the present work.

\subsection*{Acknowledgments}

We wish to thank Xavier Méla for having shown us the beautiful shape of $B_{2n,n}$, which
prompted our interest in this topic, and Jon Aaronson and G\'erard Grancher who informed us that
our curve was the well-known Blancmange function.

\section{Results}

\label{limit_standard}
As mentioned before, we have to 
renormalize $F_{n,k}$ into a new function $\varphi_{n,k}$ defined
on $[0,1]$, vanishing at 0 and 1, and vertically scaled so that the point
corresponding to the end of $B_{n-1,k-1}$ is mapped to 1:
\begin{equation}
  \label{deF_ren}
  \varphi_{n,k}(t)\ :=\ \frac{ F_{n,k}\Bigl(t\cdot\binom{n}{k}\Bigr) - tF_{n,k}\Bigl(\binom{n}{k}\Bigr) }
{ F_{n,k}\Bigl(\binom{n-1}{k-1}\Bigr) - \frac{\binom{n-1}{k-1}}{\binom{n}{k}} F_{n,k}\Bigl(\binom{n}{k}\Bigr) }\ .
\end{equation}

\begin{theorem}
\label{thm_convergence1}
  Let $\varphi_{n,k}$ be the renormalized curve associated to the basic block $B_{n,k}$ (see \eqref{deF_ren}).
  For any $p\in]0,1[$, for any sequence $(k(n))$ such that $k(n)/n\to p$, we have
  \begin{equation}
    \label{convergence1}
    \varphi_{n,k(n)}\ \converge{n\to\infty}{L^{\infty}} \macdo_p.
  \end{equation}
  Moreover, the denominator in~\eqref{deF_ren} is of order
  $ \frac{1}{n}\binom{n}{k(n)}. $
\end{theorem}

\subsection{Ergodic interpretation}
\label{ErgodicInterpretation}
The functions $F_{n,k}$ and $\varphi_{n,k}$ introduced before can be interpreted as particular cases of the following general situation: Consider a real-valued function $g$ defined on a probability space $(X,\mu)$ on which acts a measure-preserving transformation $T$. Given a point $x\in X$ and an integer $\ell\ge 1$, we construct the continuous function $F_{x,\ell}^g\,:\ [0,\ell]\to\mathbb{R}$ by $F_{x,\ell}^g(0):=0$;
for each integer $j$, $1\le j\le \ell$,
\begin{equation}
\label{defFg}
F_{x,\ell}^g(j)\ :=\ \sum_{k=0}^{j-1} g\left( T^kx\right);
\end{equation} 
and $F_{x,\ell}^g$ is linearly interpolated between the integers. 

If $g$ is integrable (which we henceforth assume), the ergodic theorem implies that, for $0<t<1$, for almost every $x$,
$$ \lim_{\ell\to\infty} \frac{1}{\ell}\sum_{0\le j < t\ell} g\left( T^jx\right)\ 
   =\ t\lim_{\ell\to\infty} \frac{1}{\ell}\sum_{0\le j < \ell} g\left( T^jx\right).
$$
Therefore, when dividing by $\ell$ both the abscissa and the ordinate, the graph of $F_{x,\ell}^g$ for large $\ell$ looks very much like a straight line with slope the empirical mean ${1}/{\ell}\sum_{0\le j < \ell} g\left( T^jx\right)$. If we want to study small fluctuations in the ergodic theorem, it is natural to remove the dominant contribution of this straight line, and rescale the ordinate to make the fluctuations appear. This leads to introduce a renormalized function $\varphi_{x,\ell}^g\,:\ [0,1]\to\mathbb{R}$ by setting 
\begin{equation*}
\label{defphig}
\varphi_{x,\ell}^g(t)\ :=\ \dfrac{F_{x,\ell}^g(t\ell)-tF_{x,\ell}^g(\ell)}{R^g_{x,\ell}},
\end{equation*} 
where $R^g_{x,\ell}$ is the renormalization in the $y$-direction, which we can canonically define by
\begin{equation}
\label{defRxl}
R^g_{x,\ell}\ :=\ \begin{cases}
             \max_{0\le t\le 1} |F_{x,\ell}^g(t\ell)-tF_{x,\ell}^g(\ell)
                              | & \text{provided this quantity does not vanish,}\\
               1         & \text{otherwise.}
              \end{cases}
\end{equation} 
It will be useful in the sequel to note that $\varphi_{x,\ell}^g$ is not changed when we add a constant to $g$.

If we consider a Bernoulli shift in which the functions $g\circ T^k$ are i.i.d. random variables, Donsker
invariance principle shows that these corrections to the law of large numbers are given by a suitably scaled
Brownian bridge. 

We are going to investigate the corresponding questions
in the context of the Pascal-adic transformation.  Let us recall, see Appendix~\ref{constructionPA}, that the sequence of letters $a$ and $b$ in the word $B_{n, k}$ encodes the trajectory $\left( x,Tx,\ldots, T^{\binom{n}{k}-1}x\right) $ of a point $x$ lying in the basis of the tower $\tau_{n,k}$ with respect to the partition
$[0, 1/2[$ (labelled by ``$a$") and $[1/2, 1[$ (labelled by ``$b$"). Thus, the function $F_{n,k}$ is nothing else but $F_{x,\binom{n}{k}}^g$ for $x$ in the basis of $\tau_{n,k}$, with the function $g$ defined by
\begin{equation*}
\label{defg}
g\ :=\ \ind{[0,1/2[}-\ind{[1/2,1[}.
\end{equation*} 
Now, the vertical renormalization chosen to define $\varphi_{n,k}$ was not exactly the one defined by~\eqref{defRxl}, but it is not difficult to restate Theorem~\ref{thm_convergence1} in the following way, where we define
\begin{equation}
\label{varphi}
\varphi_{n,k}^g\ :=\ \varphi_{x,\binom{n}{k}}^g 
\end{equation}
for any point $x$ in the  basis of $\tau_{n,k}$.

\begin{theorem}
\label{thm_convergence1bis}
Let $g\ =\ \ind{[0,1/2[}-\ind{[1/2,1[}$. If 
$k(n)/n\to p$, then
\begin{equation}
    \label{convergence1bis}
    \varphi_{n,k(n)}^g\ \converge{n\to\infty}{L^{\infty}}
    \ \dfrac{\macdo_p}{\Vert \macdo_p \Vert_\infty}. 
  \end{equation}
\end{theorem}
This result shows that the corrections to the ergodic theorem are
given by a deterministic function, when we consider sums along the 
Rokhlin towers $\tau_{n,k}$. It is possible to derive an analogous pointwise
statement at the cost of extracting a subsequence. 
\begin{theorem}
\label{thm_convergence1ter}
Let $g\ =\ \ind{[0,1/2[}-\ind{[1/2,1[}$. For $\mu_p$ almost every $x\in X$,
there exists a sequence $\ell_n$ such that 
\begin{equation}
    \label{convergence1ter}
    \varphi_{x,\ell_n}^g\ \converge{n\to\infty}{L^{\infty}}
    \ \dfrac{\macdo_p}{\Vert \macdo_p \Vert_\infty}. 
  \end{equation}
\end{theorem}

\subsection{Limit for general dyadic functions}

Let $N_0 \geq 1$. We consider the dyadic partition $\P_{N_0}$ of the interval $[0, 1[$ into $2^{N_0}$ sub-intervals. 
We want to extend Theorem~\ref{thm_convergence1bis} to a general $\P_{N_0}$-measurable real-valued function $g$. 

Associated to a $\P_{N_0}$-measurable function $g$, we define the words $B^{N_0}_{n,k}$ ($n\ge N_0,\ 0\le k\le n$) on an alphabet with $N_0+1$ letters, $\{a_0, \dots, a_{N_0}\}$. They are inductively defined by 
$B^{N_0}_{N_0,k}\ :=\ a_k $ for $0\leq k\leq N_0$, 
$B^{N_0}_{n, 0}\ := \ a_0$, $B^{N_0}_{n,n}\ := \ a_{N_0}$ for $n\geq N_0$, 
and for $0<k<n$
$$B^{N_0}_{n,k}\ :=\ B^{N_0}_{n-1,k-1}B^{N_0}_{n-1,k}. $$
To each letter $a_k$ corresponds a continuous function $F_{N_0,k}^g$ defined for $\ell$ integer in $[0,\binom{N_0}{k}]$ by 
$$ F_{N_0,k}^g(\ell)\ := \ \sum_{j=0}^{\ell-1} g(T^jx), $$
where $x$ is any point in the basis of $\tau_{N_0,k}$, and extended to the full interval by linear interpolation.

As before, we associate to the word $B^{N_0}_{n,k}$ a continuous function $F^g_{n,k} : [0, \binom{n}{k}] \to \mathbb{R}$ such that $F^g_{n, k}(0)=0$. Its graph is constructed by substituting to each letter $a_k$ the (translated) graph of $F_{N_0,k}^g$.

A central characteristic of $g$ is given by its ergodic sums
along the towers $\tau_ {N_0, k}$:
$$ {h}^g_{N_0, k}\ :=\  F_{N_0,k}^g(\tbinom{N_0}{k}), $$
for $k=0,\ldots,N_0$.

We are interested in the renormalized function $\varphi^g_{n,k}$ defined as in
\eqref{varphi} but now for a general $g$. Denoting by $R_{n,k}^g$ the renormalization
constant 
\begin{equation}
\label{defRnk}
R^g_{n,k}\ :=\ \begin{cases}
             \max_{0\le t\le 1} |F_{n,k}^g\left(t\tbinom{n}{k}\right)-tF_{n,k}^g\left(\tbinom{n}{k}\right)
                              | & \text{provided this quantity does not vanish,}\\
               1         & \text{otherwise,}
              \end{cases}
\end{equation} 
the function $\varphi^g_{n,k}$ is given, for $0\le t\le 1$, by
$$\varphi^g_{n,k}(t)\ =\ \dfrac{F_{n,k}^g\left(t\tbinom{n}{k}\right)-tF_{n,k}^g\left(\tbinom{n}{k}\right)}{R^g_{n,k}}. 
$$

As a first observation, we can point out that it is easy to find  $\P_{N_0}$-measurable functions for which convergence of $\varphi^g_{n,k}$
to a continuous function will never hold.
 
\begin{lemma}
\label{coboundary-lemma}
Let $g$ be a $\P_{N_0}$-measurable coboundary of the form $g=\h-\h\circ T$ for some bounded measurable function $\h$. Then %$h^g_{N_0, k}=0$ for all $k=0,\ldots,N_0$, and 
\begin{equation}
\label{borne}
\sup_{n,k} \|F_{n,k}^g\|_\infty\ <\ +\infty. 
\end{equation}
If $g$ is not identically $0$, then there is no cluster point in $L^\infty([0,1])$
for any sequence $\varphi_{n,k(n)}^g$ with $k(n)/n \to p\in(0,1)$. 
\end{lemma}

Note that coboundaries such as those appearing in the statement of the lemma do really exist. A simple example is given by $g:=\ind{[1/4,1/2[}-\ind{[1/2,3/4[}$, with 
transfer function $\h:=-\ind{[1/2,3/4[}$. Also note that the conclusion of the lemma 
still holds for $g$ cohomologous to a constant in $L^{\infty}$,  \emph{i.e.} of the form
$$ g\ =\ \h-\h\circ T + C $$
with $\h$ bounded measurable.
This follows from the fact that $\varphi_{n,k}^g$ is unchanged when we add a constant to $g$.

\begin{theorem}
\label{convergenceinlinfinity}
Let $g$ be measurable with respect to the dyadic partition $\P_{N_0}$. We suppose that $g$ is not cohomologous to a constant in $L^{\infty}$.  
 For any sequence $(k(n))$ such that $k(n)/n\to p\in(0,1)$, we can extract a subsequence $(n_s)$ such that  
 $\varphi^g_{n_s,k(n_s)}$ converges in $L^\infty$ to a continuous function.
\end{theorem}

In the course of the proof of this theorem, we will establish the following characterization
of $\P_{N_0}$-measurable functions $g$ which are cohomologous to some constant
in $L^\infty$:
\begin{lemma} 
 Let $g$ be $\P_{N_0}$-measurable, then $g=C+\h-\h\circ T$ for some constant $C$ and some bounded function $\h$ if and only if the quantities $h_{N_0, \ell}^g$ ($0\le \ell\le N_0$)
are proportional to $\binom{N_0}{\ell}$ ($0\le \ell\le N_0$).
\end{lemma}

Our final result concerns the cluster points we can get for the functions $\varphi_{n,k(n)}^g$. Surprisingly, the self-affine maps $\macdo_p$ which arose in the study of the basic blocks $B_{n,k}$ turn out to be the only possible limit in 
``almost all" cases for a dyadic function $g$, in a sense made clear by the following theorem. Before stating it, we introduce for any $\P_{N_0}$-measurable function $g$
the polynomial in~$p$
\begin{equation}
\label{P^g}
P^g(p)\ :=\ \sum_{\ell=0}^{N_0} h_{N_0,\ell}^g \, p^\ell(1-p)^{N_0-\ell}(N_0p-\ell).
\end{equation} 

Note that $P^g\not\equiv 0$ if and only if $g$ is not cohomologous to a constant in $L^{\infty}$. Indeed, setting $q:=p/(1-p)$ it is easy to compute the coefficients
of the corresponding polynomial in $q$ and see that they vanish if and only if 
$h_{N_0,\ell}^g\propto \binom{N_0}{\ell}$. 

\begin{theorem}
\label{thm_general}
Let $g$ be measurable with respect to the dyadic partition $\P_{N_0}$. 
%We suppose that $g$ is not cohomologous to a constant in $L^{\infty}$.  
If $P^g(p)\neq0$, for any sequence $(k(n))$ such that $k(n)/n\to p$, we have
  \begin{equation}
    \label{convergence2}
    \varphi^g_{n,k(n)}\ \converge{n\to\infty}{L^{\infty}} \sign(P^g(p))\ \macdo_p/\|\macdo_p\|_{\infty}.
  \end{equation}
Moreover, $R_{n,k(n)}^g$ is in this case of order $\frac{1}{n}\binom{n}{k(n)}$.
\end{theorem}

\subsection{Some examples where another curve appears}

Theorem~\ref{thm_general} does not characterize the possible limits for the 
values of $p$ where the polynomial $P^g$ vanishes. 
We do not have any general result in that case, but we present in Section~\ref{AnotherCurve} the study 
of some particular cases showing that some other curves can appear.

\section{Proofs}

\subsection{Piecewise linear graph associated with a triangular array}

Until now we have always considered triangular arrays of Pascal type, in which positions are denoted by couples $(n,k)$, with line $(n+1)$ traditionally
represented below line $n$. We call such arrays ``descending". We are also going to use another type of triangular arrays, which we call ``ascending", in which coordinates will be denoted by $(i,j)$, with line $(i+1)$ above line $i$. In both cases, the object located at a given position in the array is obtained from the two objects just above it by summation or concatenation.

We consider here an ascending triangular array $\ta$ with a finite number of lines, labelled (from bottom to top) $0,1,\ldots,m$.
For each $i\in\{0,\ldots,m\}$, line $i$ is constituted by $(i+1)$ pairs of real numbers $(x_{i,0},y_{i,0}),\ldots,(x_{i,i},y_{i,i})$,
satisfying the following properties:
\begin{itemize}
\item for $0\le j\le i\le m$, $x_{i,j}>0$; %, and $x_{0,0}=1$ ??????
\item for $0\le j\le i < m$, $x_{i,j}=x_{i+1,j}+x_{i+1,j+1}$ and $y_{i,j}=y_{i+1,j}+y_{i+1,j+1}$.
\end{itemize}
Observe that, because of the additive relation
existing between these numbers, we can recover the whole array if we only know its values on one of its side
(for example, if we know all the pairs $(x_{i,0},y_{i,0})$ for $0\le i\le m$). 

These pairs of real numbers are interpreted as the horizontal and vertical displacements between points on the graph
of some piecewise linear map $\varphi_\sta$. The map $\varphi_\sta$ is defined inductively on $[0,x_{0,0}]$
in the following way. First, corresponding to line 0,
we set $\varphi_\sta(0):=0$ and $\varphi_\sta(x_{0,0}):=y_{0,0}$. In general, taking into account all lines up to line $i$
provides a subdivision of $[0,x_{0,0}]$ into $2^i$ intervals whose lengths are the $x_{i,j}$ (taken several times), and
defines $\varphi_\sta$ on the bounds of these intervals. Let $I=[t,t+x_{i,j}]$ be such an interval defined in line $i$:
we must have 
$$ \varphi_\sta(t+x_{i,j})\ =\ \varphi_\sta(t) + y_{i,j}. $$
Then, coming to line $(i+1)$, $I$ is subdivised into two subintervals $[t,t+x_{i+1,j}]$ and $[t+x_{i+1,j},t+x_{i+1,j}+x_{i+1,j+1}]$
and we set 
$$\varphi_\sta(t+x_{i+1,j}):=\varphi_\sta(t) + y_{i+1,j}.$$
This inductive procedure defines at the end the values of $\varphi_\sta(t)$ for all bounds $t$ of some subdivision of $[0,x_{0,0}]$
into $2^m$ subinterval. At last, $\varphi_\sta$ is linearly interpolated between these bounds.
\begin{figure}
   \centering
   \begin{picture}(0,0)%
\includegraphics{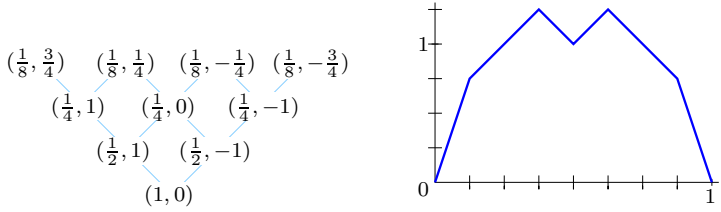}%
\end{picture}%
\setlength{\unitlength}{1906sp}%
\begingroup\makeatletter\ifx\SetFigFont\undefined%
\gdef\SetFigFont#1#2#3#4#5{%
  \reset@font\fontsize{#1}{#2pt}%
  \fontfamily{#3}\fontseries{#4}\fontshape{#5}%
  \selectfont}%
\fi\endgroup%
\begin{picture}(9505,2680)(4275,-6464)
\put(6573,-5804){\makebox(0,0)[lb]{\smash{{\SetFigFont{8}{9.6}{\rmdefault}{\mddefault}{\updefault}{\color[rgb]{0,0,0}$(\tfrac12,-1)$}%
}}}}
\put(9676,-6316){\makebox(0,0)[lb]{\smash{{\SetFigFont{8}{9.6}{\rmdefault}{\mddefault}{\updefault}{\color[rgb]{0,0,0}$0$}%
}}}}
\put(9676,-4426){\makebox(0,0)[lb]{\smash{{\SetFigFont{8}{9.6}{\rmdefault}{\mddefault}{\updefault}{\color[rgb]{0,0,0}$1$}%
}}}}
\put(13411,-6406){\makebox(0,0)[lb]{\smash{{\SetFigFont{8}{9.6}{\rmdefault}{\mddefault}{\updefault}{\color[rgb]{0,0,0}$1$}%
}}}}
\put(5491,-5804){\makebox(0,0)[lb]{\smash{{\SetFigFont{8}{9.6}{\rmdefault}{\mddefault}{\updefault}{\color[rgb]{0,0,0}$(\tfrac12,1)$}%
}}}}
\put(4906,-5233){\makebox(0,0)[lb]{\smash{{\SetFigFont{8}{9.6}{\rmdefault}{\mddefault}{\updefault}{\color[rgb]{0,0,0}$(\tfrac14,1)$}%
}}}}
\put(7201,-5233){\makebox(0,0)[lb]{\smash{{\SetFigFont{8}{9.6}{\rmdefault}{\mddefault}{\updefault}{\color[rgb]{0,0,0}$(\tfrac14,-1)$}%
}}}}
\put(4321,-4651){\makebox(0,0)[lb]{\smash{{\SetFigFont{8}{9.6}{\rmdefault}{\mddefault}{\updefault}{\color[rgb]{0,0,0}$(\tfrac18,\tfrac34)$}%
}}}}
\put(6121,-6376){\makebox(0,0)[lb]{\smash{{\SetFigFont{8}{9.6}{\rmdefault}{\mddefault}{\updefault}{\color[rgb]{0,0,0}$(1,0)$}%
}}}}
\put(6076,-5233){\makebox(0,0)[lb]{\smash{{\SetFigFont{8}{9.6}{\rmdefault}{\mddefault}{\updefault}{\color[rgb]{0,0,0}$(\tfrac14,0)$}%
}}}}
\put(5491,-4661){\makebox(0,0)[lb]{\smash{{\SetFigFont{8}{9.6}{\rmdefault}{\mddefault}{\updefault}{\color[rgb]{0,0,0}$(\tfrac18,\tfrac14)$}%
}}}}
\put(6571,-4661){\makebox(0,0)[lb]{\smash{{\SetFigFont{8}{9.6}{\rmdefault}{\mddefault}{\updefault}{\color[rgb]{0,0,0}$(\tfrac18,-\tfrac14)$}%
}}}}
\put(7786,-4661){\makebox(0,0)[lb]{\smash{{\SetFigFont{8}{9.6}{\rmdefault}{\mddefault}{\updefault}{\color[rgb]{0,0,0}$(\tfrac18,-\tfrac34)$}%
}}}}
\end{picture}%
   \caption{An array and its associated piecewise linear graph (actually this is the array $\ta_{1/2,3}$ providing the $3$rd stage approximation to $\macdo_{1/2}$).}
   \label{fig:array}
\end{figure}
\medskip

Given the triangular array $\ta$ with $m+1$ lines, we can 
construct two smaller arrays with $m$ lines denoted by $\tal$ and $\tar$: For $0\le i\le m-1$, line $i$ of
$\tal$ is constitued by the $(i+1)$ first pairs of reals in line $(i+1)$ of $\ta$, and  line $i$ of
$\tar$ is constitued by the $(i+1)$ last pairs of reals in line $(i+1)$ of $\ta$. In the sequel, we will make 
use of the following fact, whose verification is left to the reader: The graph of $\varphi_\sta$ is formed
by putting together the graphs of $\varphi_\star$ and $\varphi_\stal$. More precisely, we have
\begin{equation*}
  \label{L_and_R}
  \varphi_\sta(t)\ =\
  \begin{cases}
    \varphi_\star(t) & \text{if $0\le t\le x_{1,0}$,}\\
    \varphi_\star(x_{1,0}) + \varphi_\stal(t-x_{1,0}) & \text{if $x_{1,0}\le t\le x_{0,0}$.}
  \end{cases}
\end{equation*}

\medskip

We say that a map $\varphi\ :\ [0,x_{0,0}]\to\mathbb{R}$ is \emph{compatible} with the array $\ta$ if $\varphi(t)=\varphi_\sta(t)$
for all bound $t$ of the subdivision defined by the array.

\begin{lemma}
\label{defmacdo}  
For all $0<p<1$ and all $m\ge 0$, the self-affine map $\macdo_p$ is compatible with the triangular array $\ta_p^m$ defined
  by its lower-left side as follows: for each $0\le i\le m$,  $x_{i,0}:=p^i$ and $y_{i,0}:=ip^{i-1}$.
\end{lemma}

\begin{proof}
  We consider two transformations $\lambda_L$ and $\lambda_D$ of $\mathbb{R}^2$, which are the
respective linear parts of the affine maps $\alpha_L$ and $\alpha_R$ arising in the definition of $\macdo_p$:
$$ (x,y)\ \mathop{\longmapsto}^{\lambda_L}\ (px,py+x), $$
and
$$ (x,y)\ \mathop{\longmapsto}^{\lambda_R}\ \Bigl((1-p)x,(1-p)y-x\Bigr). $$
It is easy to check that in the triangular array $\ta_p^m$, we have for $0\le i<m$ and $0\le j\le i$
$$ (x_{i+1,j},y_{i+1,j})\ =\ \lambda_L(x_{i,j},y_{i,j}), $$
and
$$ (x_{i+1,j+1},y_{i+1,j+1})\ =\ \lambda_R(x_{i,j},y_{i,j}). $$
{From} this, we can deduce that the left part of the graph of $\varphi_{\sta_p^m}$, which is the graph of $\varphi_{\star_p^m}$,
is the image of the graph of $\varphi_{\sta_p^{m-1}}$ by the affine map $\alpha_L$, and the right part of the graph
of $\varphi_{\sta_p^m}$ is the image of the graph of $\varphi_{\sta_p^{m-1}}$ by the affine map $\alpha_R$.
A simple induction on $n$ then gives the result stated in the lemma.
\end{proof}

\subsection{Proof of Theorem \ref{thm_convergence1}}

For any $m <n $, the block $B_{n,k}$ is the concatenation of $2^m$ subblocks $B_{n-m,\cdot}$. Let us denote by $t_{n,k,m}^r$, $r = 1, \ldots, 2^m$ the position of the last letter of the $r$th subblock in the block $B_{n,k}$.
We also denote by $h_{n,k}$ the height of the basic block $B_{n,k}$, i.e. the difference between the numbers of $a$ and $b$ appearing in $B_{n,k}$. The function $\varphi_{n,k}$ is compatible with the array $\ta_{n,k}^{m}$ defined by its lower-left side as follows: For each $0\leq i \leq m$,
$$
x^{n,k}_{i,0} = t_{n,k,i}^1 \Big/\binom{n}{k} = \binom{n-i}{k-i}\Big/\binom{n}{k}\,,
$$
and
$$
y^{n,k}_{i,0} = \frac{ h_{n-i,k-i} - x^{n,k}_{i,0} h_{n,k} } { h_{n-1,k-1} - x^{n,k}_{1,0} h_{n,k} }\,.
$$

\begin{lemma}
\label{lem_ConvOfArray}
For any $m\geq 0$, any $0<p<1$ and any sequence $k(n)$ such that $\lim_n k(n)/n = p$, we have that
$$
\lim_{n\to\infty} \ta_{n,k(n)}^{m} = \ta_p^m\,,
$$
where $\ta_p^m$ was introduced in Lemma~\ref{defmacdo}.
\end{lemma}
\begin{proof}[Proof of Lemma~\ref{lem_ConvOfArray}]
It is of course sufficient to prove the convergence for the elements appearing in the lower-left side of $\ta_{n,k(n)}^{m}$. We first have
$$
\lim_{n\to\infty} x^{n,k(n)}_{i,0} = \lim_{n\to\infty} \prod_{r=0}^{i-1}\frac{k(n)-r}{n-r} = p^i\,.
$$
Moreover, using the identity $h_{n,k} = \frac{n-2k}{n} \binom{n}{k}$, we also obtain
\begin{align*}
\lim_{n\to\infty} y^{n,k(n)}_{i,0}
&= \lim_{n\to\infty} \frac{ h_{n-i,k(n)-i} - x^{n,k}_{i,0} h_{n,k(n)} } { h_{n-1,k(n)-1} - x^{n,k}_{1,0} h_{n,k(n)} }\\
&= \lim_{n\to\infty} \frac{\binom{n-i}{k(n)-i} \left( \frac{n+i-2k(n)}{n-i} - \frac{n-2k(n)}{n} \right)}{\binom{n-1}{k(n)-1} \left( \frac{n+1-2k(n)}{n-1} - \frac{n-2k(n)}{n} \right)}\\
&= \lim_{n\to\infty} i\, \frac{n-1}{n-i} \prod_{r=1}^{i-1} \frac{k(n)-r}{n-r}\\
&= i p^{i-1}\,.
\end{align*}
\end{proof}

We denote by $\varphi_{n,k}^{m} := \varphi_{\sta_{n,k}^{m}}$ the polygonal approximation of $\varphi_{n,k}$ at the order $m$.
A computation shows that $\varphi_{n,k} = \varphi_{n,k}^{n-1}$.

Lemma~\ref{lem_ConvOfArray} obviously implies
\begin{equation}
\lim_{m\to +\infty}\limsup_{n\to +\infty}\Vert\varphi_{n,k(n)}^{m}-\varphi_{\sta_{p}^m}\Vert_{\infty} = 0.
\end{equation}
Moreover, by continuity of $\macdo_p$, we also have
\begin{equation}
\varphi_{\sta_{p}^m}\ \converge{m\to\infty}{L^{\infty}} \macdo_p.
\end{equation}
Hence, it is enough to prove that
\begin{equation}
\lim_{m\to +\infty}\limsup_{n\to +\infty}\Vert\varphi_{n,k(n)}^{m}-\varphi_{n,k(n)}\Vert_{\infty} = 0,
\end{equation}
which is a consequence of the general Theorem~\ref{convergenceinlinfinity}.
\qed

\subsection{Proof of Theorem \ref{thm_convergence1ter}}

For each $n\ge1$, let us denote by $k_n(x)$ the unique index such that
$$ x\ \in\ \tau_{n,k_n(x)}. $$
We have
$$ \dfrac{k_n(x)}{n}\ \converge{n\rightarrow\infty}{}\ p\qquad \text{$\mu_p$-almost surely}.
$$
Thus, it follows from Theorem~\ref{thm_convergence1bis} that
$$
\varphi_{n,k_n(x)}^g\ \converge{n\to\infty}{L^{\infty}}
    \ \dfrac{\macdo_p}{\Vert \macdo_p \Vert_\infty} \qquad \text{$\mu_p$-almost surely}.
$$
It therefore only remains to observe that $\mu_p$-almost surely, $x$ lies 
arbitrarily close to the bottom of $\tau_{n,k_n(x)}$ for infinitely many $n$;
more precisely there exists a sequence $(n_s)$ such that the height of $x$ in tower
$\tau_{n_s,k_{n_s}(x)}$ is smaller than $\frac{1}{s}\binom{n_s}{k_{n_s}(x)}$. This follows from \cite[Lemma 2.5]{Janvresse-delaRue04}.
\qed
 
\subsection{Proof of Lemma~\ref{coboundary-lemma}}

For any $\ell\in\{0,\ldots,\binom{n}{k}\}$, 
$$ F_{n,k}^g(\ell)\ =\ \h(x)-\h(T^{\ell}x) $$
for any $x$ in the basis of the tower $\tau_{n,k}$, which proves~\eqref{borne}.
This implies that the renormalization constants $R_{n,k}^g$ are uniformly bounded. 
It is easy to see that for $n$ large enough, each letter $a_k$ appears at least once 
in the decomposition of the word $B_{n,k(n)}^{N_0}$ (here we use the assumption that $\lim k(n)/n\in(0,1)$). 
Suppose first that there exists $0\le k_0\le N_0$ such that $g$ is not constant on $\tau_{N_0,k_0}$. Then on each subinterval of $[0,1]$ corresponding to one occurence
of $a_k$ in $B_{n,k(n)}^{N_0}$, the function $\varphi_{n,k(n)}^g$ has {variation}
which is uniformly bounded below by some $c>0$. The conclusion follows since the length of this subinterval goes to 0 as $n\to \infty$. 

Finally suppose that $g$ is constant on each tower $\tau_{N_0,k}$. Note that this 
constant cannot be the same for every towers otherwise $g$ would be identically 0
(remember that $g$ is a coboundary). Hence there exists $k_1$ such that $g$ takes different values on $\tau_{N_0,k_1}$ and $\tau_{N_0,k_1+1}$. Therefore $g$ is not constant on the tower $\tau_{N_0+1,k_1}$ and we are back to the previous case.
\qed

\subsection{Proof of Theorem \ref{convergenceinlinfinity}}

The function $\varphi^g_{n,k(n)}%=\varphi_{\sta_{n, k(n)}^g}
$ is compatible with the array $\ta_{n, k(n)}^g$, which is defined by its lower-left side: for $0\le i\le n-N_0$,
\begin{equation}
\label{def_x^gy^g}
x_{i, 0}^{n, k(n), g} := \dfrac{\binom{n-i}{k(n)-i}}{\binom{n}{k(n)}}, \qquad 
y_{i, 0}^{n, k(n), g}  := \varphi^g_{n,k(n)}\left(x_{i, 0}^{n, k(n), g}\right). 
\end{equation}
Moreover, 
\begin{equation}
\label{Rgrand}
\|\varphi^g_{n,k(n)}-\varphi_{\sta_{n, k(n)}^g}\|_{\infty}\ \converge{n\to\infty}{}\ 0,
\end{equation} 
provided that the renormalization constant $R_{n,k(n)}^g$ goes to infinity. (This means that in this case we can forget the variations in each $F_{N_0,k}^g$ and replace them by linear functions.)

Notice that the $y$ coefficients %$y_{n-N_0,j}^{n, k(n), g}$ 
on the top line of the ascending array $\ta_{n, k(n)}^g$ are either null or of the form $\alpha_\ell(n,k(n))/R_{n, k(n)}^g$, for $0\le \ell\le N_0$, where 
$$
\alpha_\ell(n,k(n))\ =\  h_{N_0, \ell}^g - \binom{N_0}{\ell} \sum_{r=0}^{N_0} h_{N_0, r}^g \dfrac{\binom{n-N_0}{k(n)-r}}{\binom{n}{k(n)}}.
$$
The quantity substracted to $h_{N_0, \ell}^g$ corresponds to adding a constant $d$
to $g$ so that $F_{n,k(n)}^{g+d}$ vanishes at its end point. Thus, we can rewrite 
$y_{i, j}^{n, k(n), g}$ as 
$$
y_{i, j}^{n, k(n), g}  = 
\dfrac{1}{R_{n, k(n)}^g}\sum_{\ell=0}^{N_0}\alpha_\ell(n,k(n))\binom{n-i-N_0}{k(n)-i-\ell+j}.
$$
We denote by $\ta_{n, k(n)}^{g, m}$ the truncated array constituted by the first $(m+1)$ lines of $\ta_{n, k(n)}^{g}$. 
Observe that the coefficients $y_{i, j}^{n, k(n), g}$ satisfy the conditions of Proposition~\ref{prop} stated below, with $\delta:=\min\{p/4,(1-p)/4\}$. 
In particular, it follows from the latter that 
$$
\sup_{0\le j\le i} |y_{i, j}^{n, k(n), g}| \le 3 e^{-Ci}, 
$$
provided that $2\delta < k(n)/n < 1-2\delta$, which is true for $n$ large enough since $p\in(0,1)$. 
This implies that 
\begin{equation*}
\label{P1}
\sup_n\| \varphi_{\sta_{n, k(n)}^g} - \varphi_{\sta_{n, k(n)}^{g, m}} \|_{\infty} \le 3\sum_{i\ge m} e^{-Ci}
\end{equation*}
which goes to zero as $m$ goes to infinity. 
For any fixed $m$, we can extract a subsequence $(n_s)$ such that the arrays $\ta_{n_s, k(n_s)}^{g,m}$ converge to an array $\ta^{g,m}$; by a classical diagonalization argument, it is then possible to find $(n_s)$ such that the convergence holds for any $m$. Equivalently, the function $\varphi_{\ta_{n_s, k(n_s)}^{g,m}}$ converges in $L^\infty$ to a function $\varphi^{g,m}$. Moreover, it follows from  Proposition~\ref{prop} that $\varphi^{g,m}$ converges in $L^\infty$ as $m$ goes to infinity to a continuous function $\varphi^{g}$. Then, provided that \eqref{Rgrand} is satisfied, we get that 
$$ \|\varphi_{n_s,k(n_s)}^g-\varphi^{g}\|_{\infty}\ \converge{n\to\infty}{}\ 0. $$
To complete the proof of the theorem, it only remains to show that if $R_{n,k(n)}^g$ 
is bounded, then $g=C+\h-\h\circ T$ with $\h$ bounded measurable.  Proposition~\ref{prop} clearly implies that $\alpha_\ell(n,k(n))/R_{n,k(n)}^g\to 0$ as $n\to\infty$.
If we assume that $R_{n,k(n)}^g$ is bounded, this leads to 
$$
h_{N_0, \ell}^g - \binom{N_0}{\ell} \gamma_n\ \converge{n\to\infty}{}\ 0
$$
for all $0\le \ell\le N_0$, where 
$$ \gamma_n\ :=\ \sum_{r=0}^{N_0} h_{N_0, r}^g \dfrac{\binom{n-N_0}{k(n)-r}}{\binom{n}{k(n)}}.
$$
This in turn implies that the quantities $h_{N_0, \ell}^g$ ($0\le \ell\le N_0$)
are proportional to $\binom{N_0}{\ell}$ ($0\le \ell\le N_0$), which means that we can
substract some constant $C$ to the function $g$ so that 
$$
h_{N_0, \ell}^{g-C}\ =\ 0\qquad \forall \ell\in\{0,\ldots,N_0\}.
$$
But this is easily seen to be equivalent to 
the following: The function $g-C$ belongs to the linear space spanned by the functions $\h_{\ell,r}-\h_{\ell,r}\circ T$ ($0\le \ell\le N_0$, $1\le r\le \tbinom{N_0}{\ell}$)
where $\h_{\ell,r}$ is the indicator function of the $r$-th rung in tower $\tau_{N_0,\ell}$.
\qed

\begin{proposition}
\label{prop}
Let $N_0\ge1$, $\delta\in(0,1/4)$, and $\nb$, $\kb$ such that $\nb\ge N_0$ and $2\delta\nb\le \kb\le(1-2\delta)\nb$.
Let $\alpha_\ell$, $\ell=0,\ldots,N_0$ be real numbers. 
For $N_0\le n\le \nb$ and $0\le k\le \kb$ with $n-k\le\nb-\kb$, we define
$$ \gamma_{n,k}\ :=\ \dfrac{1}{R}\sum_{\ell=0}^{N_0}\alpha_\ell\binom{n-N_0}{k-\ell}, $$
where $R$ is a renormalization constant such that $|\gamma_{n,k}|$ is always bounded by 2. There exists a constant $C=C(\delta,N_0)$ such that, provided $\nb$ is large enough, the following inequality holds for all $n,k$:
\begin{equation}
\label{gamma}
\gamma_{n,k}\ \le\ 3e^{-C(\nb-n)}.
\end{equation} 
\end{proposition}

\begin{proof}[Proof of Proposition~\ref{prop}]
We choose $\ell_0\in\{0,\ldots,N_0\}$. We can write 
\begin{multline*}
R\gamma_{n,k}\ =\ 
\binom{n-N_0}{k-\ell_0} \sum_{\ell=0}^{\ell_0-1} \alpha_\ell 
\prod_{\ell+1\le r\le \ell_0} \dfrac{n-N_0-k+r}{k+1-r}  \\
 + \binom{n-N_0}{k-\ell_0} \sum_{\ell=\ell_0}^{N_0} \alpha_\ell 
\prod_{\ell_0+1\le r\le \ell} \dfrac{k+1-r}{n-N_0-k+r},
\end{multline*}
provided that $\binom{n-N_0}{k-\ell_0}\neq 0$. We are going to bound the second term of the RHS; the first one can be treated in a similar way. It can be written as
\begin{equation}
\label{rhs}
\binom{n-N_0}{k-\ell_0}\ \dfrac{\tilde P(n,k)}{\tilde Q(n-k)}, 
\end{equation}
where 
$$ \tilde P(n,k)\ :=\ \sum_{\ell=\ell_0}^{N_0} \alpha_\ell
    \prod_{\ell_0+1\le r\le \ell} {(k+1-j)}
    \prod_{\ell+1\le r\le N_0} (n-N_0-k+r)
   , $$
and
$$ \tilde Q(n-k)\ :=\ 
   \prod_{\ell_0+1\le r\le N_0} (n-N_0-k+r). $$
It is convenient to make the following change of variables:
$$ x\ :=\ \kb-k\quad;\quad y\ :=\ (\nb-\kb) - (n-k), $$
and to set
$$ P(x,y)\ :=\ \tilde P(n,k)\quad ;\quad Q(y)\ :=\ \tilde Q(n-k). $$
Notice that in the domain where the $\gamma_{n,k}$ are defined, $x$ and $y$
are nonnegative integers. We observe that the degree of $P$ is $N_0-\ell_0\le N_0$, so that we can write
$$ P(x,y)\ =\ \sum_{\substack{u,v\ge0\\ u+v\le N_0}} c_{u,v}x^uy^v. $$
There exists a constant $M=M(N_0)$ such that for each polynomial $P$ of the above form, we have 
$$ \sum_{\substack{u,v\ge0\\ u+v\le N_0}} |c_{u,v}|\ \le\ M\max_{\substack{u,v\ge0\\ u+v\le N_0}}|P(u,v)|. $$
Indeed, the map $(c_{i,j})\longmapsto 
(P(u,v))$ is linear and one-to-one in a finite-dimensional space where all the norms are equivalent. Therefore, for nonnegative $x,y$
$$ |P(x,y)|\ \le\ M\max_{\substack{u,v\ge0\\ u+v\le N_0}}|P(u,v)| \ (1+x+y)^{N_0}. $$
Hence, since $Q(0)\ge Q(v)$ for any $j\ge 0$, 
$$ \left| \dfrac{P(x,y)}{Q(y)} \right|\ \le\ M\max_{\substack{u,v\ge0\\ u+v\le N_0}}\left|\dfrac{P(u,v)}{Q(v)}\right| (1+x+y)^{N_0}
\ \dfrac{Q(0)}{Q(y)}.$$
For all $y\ge 0$, we have
$$\dfrac{Q(0)}{Q(y)}\ = \prod_{r=\ell_0+1}^{N_0} \dfrac{\nb-\kb-N_0+r}{\nb-\kb-N_0-y+r}
 \le\ (1+3y)^{N_0}. $$
The last inequality is easily obtained by considering the two cases: $y<(\nb-\kb)/2$, and $y\ge(\nb-\kb)/2$. 
Therefore, we get
\begin{equation}
\label{PQ} \left|\dfrac{P(x,y)}{Q(y)}\right| \ \le\ M(1+x+y)^{N_0}(1+3y)^{N_0}\max_{\substack{u,v\ge0\\ u+v\le N_0}}\left|\dfrac{P(u,v)}{Q(v)}\right|. 
\end{equation}
We now need to estimate the maximum in the above formula. Denoting by $n_1,k_1$ the position where this maximum is attained, it follows from the assumption $\gamma_{n,k}\le 2$ that
\begin{multline}
2\ \ge\ \gamma_{n_1,k_1}\ =\ \dfrac{1}{R}\binom{n_1-N_0}{k_1-\ell_0} \max_{\substack{u,v\ge0\\ u+v\le N_0}}\left|\dfrac{P(u,v)}{Q(v)}\right| \\
\ge\ \dfrac{1}{R}\binom{\nb-N_0}{\kb-\ell_0} (2\delta)^{N_0} \max_{\substack{u,v\ge0\\ u+v\le N_0}}\left|\dfrac{P(u,v)}{Q(v)}\right|.\label{max}
\end{multline}
The last inequality follows from the fact that $k_1/n_1\in (2\delta, 1-2\delta)$. 
Therefore, provided that $\binom{n-N_0}{k-\ell_0}\neq 0$, we get from \eqref{rhs}, \eqref{PQ} and \eqref{max}
$$ \gamma_{n,k}\ \le\ C(\delta,N_0) (1+\nb-n)^{2N_0}\dfrac{\binom{n-N_0}{k-\ell_0}}{\binom{\nb-N_0}{\kb-\ell_0}}.
$$
Notice that if $(n-N_0, k-\ell_0)$ is such that $(k-\ell_0)/(n-N_0)\in (\delta,1-\delta)$, $(n-N_0, k-\ell_0)$ and $(\nb-N_0, \kb-\ell_0)$ can always be linked by a path in the triangle such that all the points along the path stay in the same set. Therefore, the result simply follows from repetition of the inequalities, valid for
$k/n \in (\delta,1-\delta)$,
\begin{gather*}
\binom{n-1}{k-1} = \frac{k}{n} \binom{n}{k}   \leq (1-\delta) \binom{n}{k}\,,\\
\binom{n-1}{k}   = \frac{n-k}{n} \binom{n}{k} \leq (1-\delta) \binom{n}{k}\,.
\end{gather*}
Suppose now that $(k-\ell_0)/(n-N_0) \le \delta$. We want to link the points $(n-N_0, k-\ell_0)$ and $(\nb-N_0, \kb-\ell_0)$ by a path staying as much as possible in the set $\{(n, k): k/n\in (\delta,1-\delta)\}$. One can easily check that the fraction of the length of such a path spent inside this set is bounded below by $1/2(1-\delta)$. Moreover, since for any $(n, k)$, $\max(\binom{n-1}{k-1}, \binom{n-1}{k})\le \binom{n}{k}$, we can repeat our argument and we obtain that 
$$
\dfrac{\binom{n-N_0}{k-\ell_0}}{\binom{\nb-N_0}{\kb-\ell_0}} \le e^{-C(\delta)(\nb-n)}.
$$
This proves our claim provided $\binom{n-N_0}{k-\ell_0}\neq 0$. However, for any $(n, k)$ it is possible to choose $0\le \ell_0\le N_0$ such that this holds. The conclusion follows since our estimate is uniform in $\ell_0$.  
\begin{figure}
   \centering
\begin{picture}(0,0)%
\includegraphics{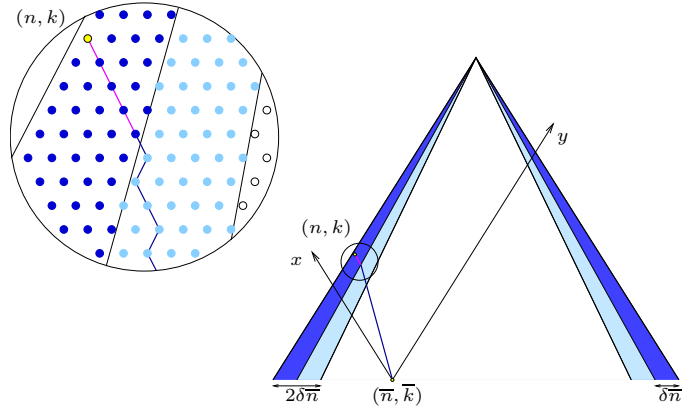}%
\end{picture}%
\setlength{\unitlength}{987sp}%
\begingroup\makeatletter\ifx\SetFigFont\undefined%
\gdef\SetFigFont#1#2#3#4#5{%
  \reset@font\fontsize{#1}{#2pt}%
  \fontfamily{#3}\fontseries{#4}\fontshape{#5}%
  \selectfont}%
\fi\endgroup%
\begin{picture}(16820,10119)(-5107,-9694)
\put(11176,-9586){\makebox(0,0)[lb]{\smash{{\SetFigFont{7}{8.4}{\rmdefault}{\mddefault}{\updefault}{\color[rgb]{0,0,0}$\delta\nb$}%
}}}}
\put(1801,-9586){\makebox(0,0)[lb]{\smash{{\SetFigFont{7}{8.4}{\rmdefault}{\mddefault}{\updefault}{\color[rgb]{0,0,0}$2\delta\nb$}%
}}}}
\put(3976,-9586){\makebox(0,0)[lb]{\smash{{\SetFigFont{7}{8.4}{\rmdefault}{\mddefault}{\updefault}{\color[rgb]{0,0,0}$(\nb,\kb)$}%
}}}}
\put(8626,-3061){\makebox(0,0)[lb]{\smash{{\SetFigFont{7}{8.4}{\rmdefault}{\mddefault}{\updefault}{\color[rgb]{0,0,0}$y$}%
}}}}
\put(1936,-6181){\makebox(0,0)[lb]{\smash{{\SetFigFont{7}{8.4}{\rmdefault}{\mddefault}{\updefault}{\color[rgb]{0,0,0}$x$}%
}}}}
\put(2191,-5371){\makebox(0,0)[lb]{\smash{{\SetFigFont{7}{8.4}{\rmdefault}{\mddefault}{\updefault}{\color[rgb]{0,0,0}$(n,k)$}%
}}}}
\put(-4971,-89){\makebox(0,0)[lb]{\smash{{\SetFigFont{7}{8.4}{\rmdefault}{\mddefault}{\updefault}{\color[rgb]{0,0,0}$(n,k)$}%
}}}}
\end{picture}%
   \caption{In the case where $\frac{k}{n}\not\in(\delta,1-\delta)$, we construct a
   path escaping as fast as possible from this region.}
   \label{fig:linfini}
\end{figure}
\end{proof}

\subsection{Proof of Theorem \ref{thm_general}}

Suppose for the beginning that, for some $0\le \ell\le N_0$, $h_{N_0,k}^g=\delta_{k,\ell}$.
Then $y_{i,0}^{n,k(n),g}$ (see \eqref{def_x^gy^g}) can be written as (writing simply $k$ for $k(n)$)
\begin{eqnarray*}
y_{i,0}^{n,k,g} & = & \dfrac{1}{R_{n,k}^g}\left(h_{n-i,k-i}^{N_0,\ell}-x_{i,0}^{n,k}h_{n,k}^{N_0,\ell}\right) \\
& = & \dfrac{1}{R_{n,k}^g}\left(\binom{n-N_0-i}{k-\ell-i}-\dfrac{\binom{n-i}{k-i}}{\binom{n}{k}}
\binom{n-N_0}{k-\ell}\right).
\end{eqnarray*}
After some algebra, we see that the numerator is equal to
\begin{eqnarray}
\label{eq_yR}
y_{i,0}^{n,k,g}R_{n,k}^g  & = & 
\binom{n-N_0}{k-\ell} \prod_{j=0}^{i-1}\dfrac{k-j}{n-j}\left(\prod_{j=0}^{i-1} \dfrac{k-\ell-j}{k-j}\dfrac{n-j}{n-N_0-j}-1\right)\\
\nonumber
& = & 
\binom{n-N_0}{k-\ell} \prod_{j=0}^{i-1}\dfrac{k-j}{n-j}\left(i\left(\frac{N_0}{n}-\frac{\ell}{k}\right)+o\left(\frac{1}{n}\right)\right) \\
\nonumber
& = & 
\binom{n-N_0}{k-\ell} \prod_{j=0}^{i-1}\dfrac{k-j}{n-j}\left(\frac{i}{k}(N_0p-\ell)+o\left(\frac{1}{n}\right)\right)\\
\nonumber
& = &  \frac{1}{n}\binom{n}{k}ip^{i-1}p^\ell(1-p)^{N_0-\ell}(N_0p-\ell+o(1)).
\end{eqnarray}
We now turn to the general case. By linearity, we get
$$ y_{i,0}^{n,k,g}R_{n,k}^g\ =\  \frac{1}{n}\binom{n}{k}ip^{i-1}P^g(p) (1+o(1)). $$
Provided that $P^g(p)\neq 0$, the denominator $R_{n,k}^g$ is proportional to the same expression where $i=1$. It follows that for some $C\neq 0$,
$$ y_{i,0}^{n,k,g}\ =\ ip^{i-1}(C+o(1)). $$
\qed

\section{Open problems and conjectures}

\subsection{Limiting curves in the transition regime}

\label{AnotherCurve}
\begin{figure}[t!]
   \centering
   \includegraphics[width=10cm,height=7cm]{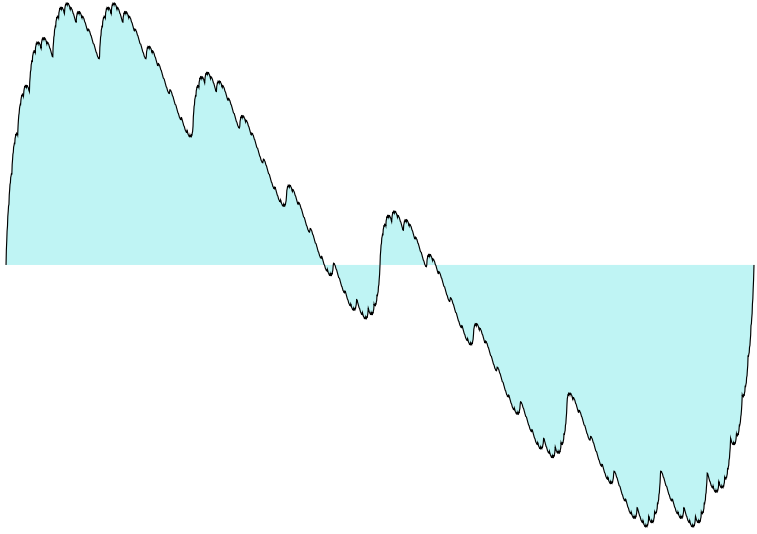}
   \caption{The limiting curve corresponding to the array $\ta_{1/2}'^{\,m}$, obtained along the sequence $(2k,k)$. Notice that this curve does not have the same self-similarity as $\macdo_p$, and is thus much less stable. For example, along the sequence $(2k-1,k-1)$ the limiting curve is the left half, which is different.} 
   \label{nacdo}
\end{figure}
Our aim here is to sudy, in some particular cases, the behaviour in the transition regime, i.e. when the polynomial $P^g$ vanishes. We introduce a family of $\P_{N_0}$-measurable functions $g_{N_0}$, indexed by $N_0=1,2,\ldots$, such that
$$
h_{N_0,\ell}^{g_{N_0}} = (-1)^\ell \binom{N_0}{\ell}\,. 
$$
It is easy to check that, for $N_0\geq 2$, $P^{g_{N_0}}$ has a zero of multiplicity $N_0-1$ at $1/2$. Indeed, using the identity $\binom{n}{k}=\binom{n-1}{k}+\binom{n-1}{k-1}$, one easily obtains that, for $N_0\geq 2$,
$$
P^{g_{N_0}}(p) = \left( 1-\frac{p}{1-p} \right) P^{g_{N_0-1}}(p) + \left( 1-2p \right)^{N_0-1}\,, 
$$
and the claim follows since $P^{g_{1}}(p) = 2p(1-p)$ does not vanish at $p=1/2$.

Using this family of functions, it is possible to investigate the behaviour of the limiting graph in the transition regime, i.e. along the sequence $(n,k)=(2k, k)$. In particular, we would like to see whether the multiplicity of the zero of the polynomial at $p=1/2$ has an influence on the limit. It turns out that, seemingly, only the parity of the multiplicity plays a crucial role.

Indeed, introducing the notation $h_{n,k}^{g_{N_0}} = S_{N_0}(n,k) \binom{n}{k}$, and using the identity $h_{n,k}^{g_{N_0}} = h_{n-1,k}^{g_{N_0-1}} - h_{n-1,k-1}^{g_{N_0-1}}$, we easily obtain the following recurrence relation,
$$
S_{N_0}(n,k) = \left( 1-\frac{k}{n} \right) S_{N_0-1}(n-1,k) - \frac{k}{n} S_{N_0-1}(n-1,k-1)\,.
$$
 From this, we can then easily compute the following asymptotics for the numerator of $y_{i,0}^{2k,k,g_{N_0}}$: It is equal to $(\frac{1}{2})^i \binom{2k}{k}$ times 
\begin{align*}
N_0 = 2 & :\quad  \phantom{-}\frac{1}{4k^2}\,i(i-1) + o(k^{-2})\,,\\
N_0 = 3 & :\quad -\frac{3}{4k^2}\,i + o(k^{-2})\,,\\
N_0 = 4 & :\quad -\frac{3}{4k^3}\,i(i-1) + o(k^{-3})\,,\\
N_0 = 5 & :\quad  \phantom{-}\frac{15}{k^3}\,i + o(k^{-3})\,,\\
N_0 = 6 & :\quad  \phantom{-}\frac{45}{16k^4}\,i(i-1) + o(k^{-4})\,,\\
N_0 = 7 & : \quad-\frac{105}{16k^4}\,i + o(k^{-4})\,,\\
\end{align*}
and so on, and so forth. We therefore see that for odd $N_0$, there is convergence to the curve $\macdo_{1/2}$, with alternating signs. Of course, the scaling is different from what we saw in the generic case. More interestingly, we see that when $N_0$ is even, there is convergence to a new curve (again with alternating signs), characterized by the array $\ta_{1/2}'^{\,m}$ defined by its lower-left side as follows: For $0\le i\le m$, $x_{i,0}':=2^{-i}$ and $y_{i,0}':=i(i-1)(1/2)^{i-2}$. A picture of this limiting curve is given in Fig.~\ref{nacdo}.

A heuristic interpretation of the preceding result is that, for the function $g_{N_0}$ that we consider here, the polynomial $P^{g_{N_0}}(p)$ changes signs when crossing $p=1/2$ for even values of $N_0$, so that we are looking at the transition between the two opposite curves $\pm \macdo_{1/2}/\| \macdo_{1/2} \|_\infty$. Actually, this can be seen in the array $\ta_{1/2}'$. Indeed, if we consider the subarray of  $\ta_{1/2}'$ starting from position $(i_0,0)$ (it is defined by its lower-left side by $x_{i,0}'^{(i_0)}:=x'_{i+i_0,0}$ and $y_{i,0}'^{(i_0)}:=y'_{i+i_0,0}$),
and renormalize the associated curve in the standard way, a little computation shows that it converges in $L^\infty$ as $i_0\to\infty$ to $\macdo_{1/2}/\| \macdo_{1/2} \|_\infty$. Proceeding similarly on the right side gives rise to the curve $-\macdo_{1/2}/\| \macdo_{1/2} \|_\infty$.

\subsubsection*{Question}
The preceding analysis leads us to raise the following question:
Is it possible to observe limiting curves other than $\pm\macdo_p$,
and portions of the curve in Fig.~\ref{nacdo}?
It is actually possible to get other curves for arbitrarily large $n$: 
for any $s\ge 1$, by taking appropriate initial conditions, 
we can get arbitrarily close to the curve given
by the array %$\ta_{p}^{(s)}$, 
defined by its lower-left side as follows: For $0\le i\le m$, $x_{i,0}^{(s)}:=p^{i}$ and $y_{i,0}^{(s)}:=i(i-1)\cdots(i-s)p^{i-s-1}$.
However, such curves do not seem to survive in the limit. 

\subsection{Larger classes of functions}

Recall that we introduced for any $\P_{N_0}$-measurable function $g$ 
the polynomial in~$p$
$$
%\label{P^g}
P^g(p)\ :=\ \sum_{\ell=0}^{N_0} h_{N_0,\ell}^g \, p^\ell(1-p)^{N_0-\ell}(N_0p-\ell),
$$
where $h_{N_0,\ell}^g $ is the sum of the values taken by $g$ on each rung of the tower $\tau_{N_0,\ell}$. 
Since $\mu_p$ gives the mass $p^\ell(1-p)^{N_0-\ell}$ to each rung of $\tau_{N_0,\ell}$, we can rewrite $h_{N_0,\ell}^g $ as 
$$
h_{N_0,\ell}^g = \dfrac{1}{p^\ell(1-p)^{N_0-\ell}}\int_{\tau_{N_0,\ell}} g d\mu_p.
$$
Therefore, the polynomial $P^g(p)$ works out to 
$$
P^g(p)\ =\ \sum_{\ell=0}^{N_0} (N_0p - \ell ) \int_{\tau_{N_0,\ell}} g(x) d\mu_p(x). 
$$
For $x\in\tau_{N_0,\ell}$, $\ell$ is equal to the sum of the first $N_0$ digits $X_1, \dots, X_{N_0}$ in the binary expansion of $x$. This allows us to rewrite the last expression as 
\begin{equation}
\label{cov}
P^g(p)
\ =\ \sum_{\ell=0}^{N_0} \int_{\tau_{N_0,\ell}} g(x) \bigl(N_0p - \sum_{i=1}^{N_0}X_i \bigr) d\mu_p(x) 
\ =\ -\ \mathrm{cov}_{\mu_p}\Bigl( g\ ;\ \sum_{i=1}^{N_0}X_i \Bigr).
\end{equation}
As we are now interested in functions which are not necessarily $\P_{N_0}$-measurable, it is convenient to emphasize the $N_0$-dependence of $P^g$ by writing $P^g_{N_0}$. 
Any $\P_{N_0}$-measurable function $g$ can also be viewed as a $\P_{N_0+1}$-measurable function. Thus, for such a function, we see from \eqref{cov} that $P^g_{N_0} = P^g_{N_0+1} $. 

For an arbitrary $g$, a natural question is the following: Suppose that 
$$\lim_{N_0\to\infty}\mathrm{cov}_{\mu_p}\Bigl( g\ ;\ \sum_{i=1}^{N_0}X_i \Bigr)$$ 
exists and is nonzero. Does the conclusion of Theorem~\ref{thm_general} still hold? 

A sufficient condition for the existence of the limit is that 
$$
\sum_{N_0}\ \Bigl\| E_{\mu_p}[g | \P_{N_0+1}] - E_{\mu_p}[g | \P_{N_0}] \Bigr\|_2 < \infty, 
$$
which is satisfied for example by indicators of intervals. 

\bigskip
It is easy to see that Theorem~\ref{thm_general} cannot hold for arbitrary measurable function $g$. Actually it
is known that in any aperiodic ergodic dynamical system, one can find a function $g$ such that the invariance principle holds \cite{volny99}; for such a function, we clearly cannot have the type of behavior described in the 
present paper. We can also construct an explicit counterexample. Let's start from $g = \ind{[0,1/2[}$, which satisfies Theorem~\ref{thm_general}. 
To each tower $\tau_{n, k}$ of a sufficiently large level $n$, we make the following procedure: We modify the values taken by $g$ at the bottom and the top of the tower. On the first $\epsilon\binom{n}{k}$ rungs, the value is set to $1$, while it is set to $0$ on the last $\epsilon\binom{n}{k}$ rungs. We repeat this construction for a sequence $\epsilon_i$ with $\sum_i\epsilon_i$ small, and levels $n_i$ chosen such that $1/n_i$ is much smaller than $\epsilon_i$. 
Since the fluctuations giving rise to $\macdo_p$ for the original function $g$ are of order $\binom{n}{k}/n$ (see Theorem~\ref{thm_general}), there cannot be convergence to $\macdo_p$ for the modified function.

\subsection{Other transformations}

\subsubsection{Generalized Pascal-adic transformations}

\begin{figure}
   \centering
   \includegraphics[height=3cm,width=4cm]{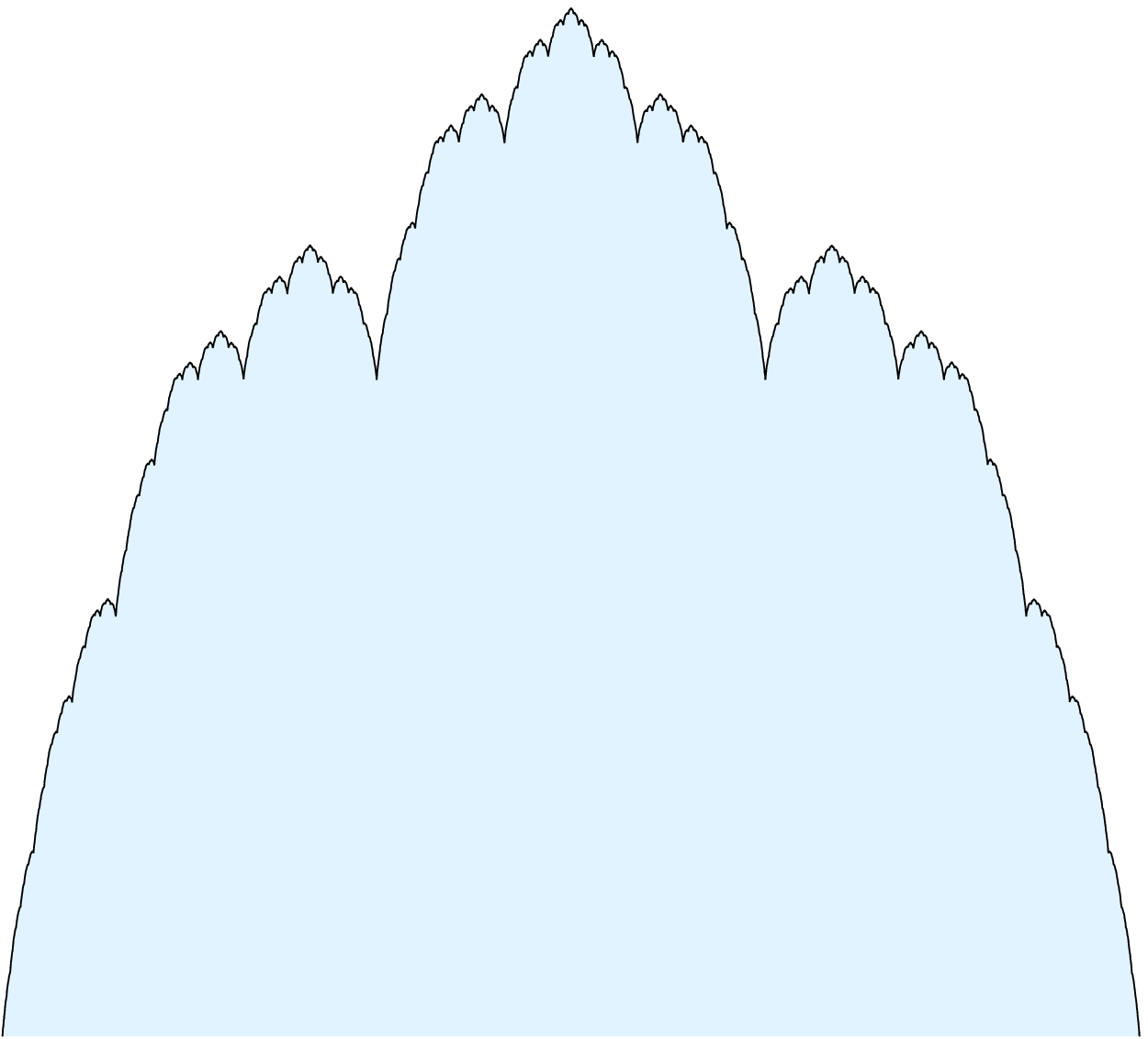}
   \hspace{2mm}
   \includegraphics[height=3cm,width=4cm]{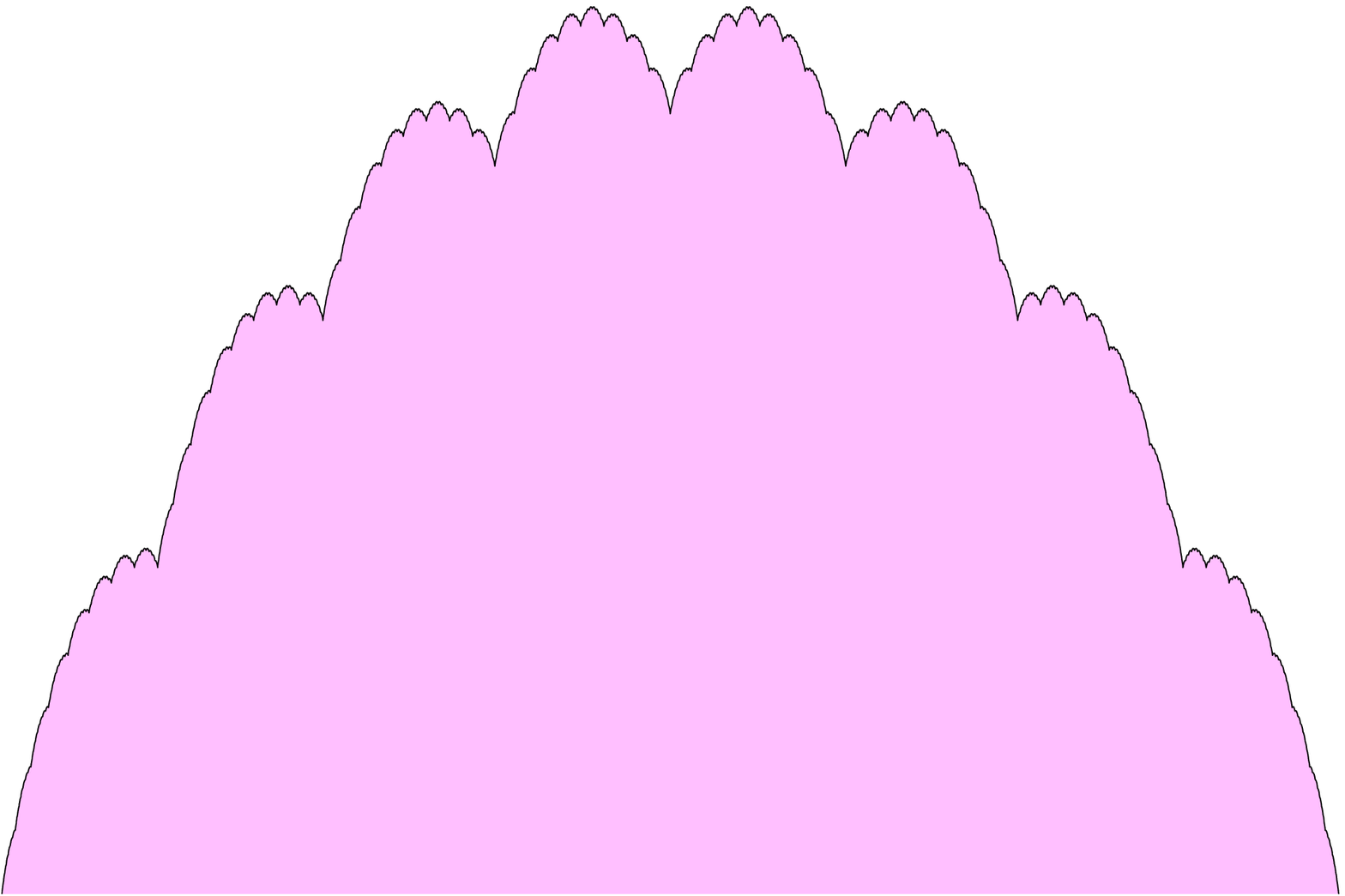}
   \hspace{2mm}
   \includegraphics[height=3cm,width=4cm]{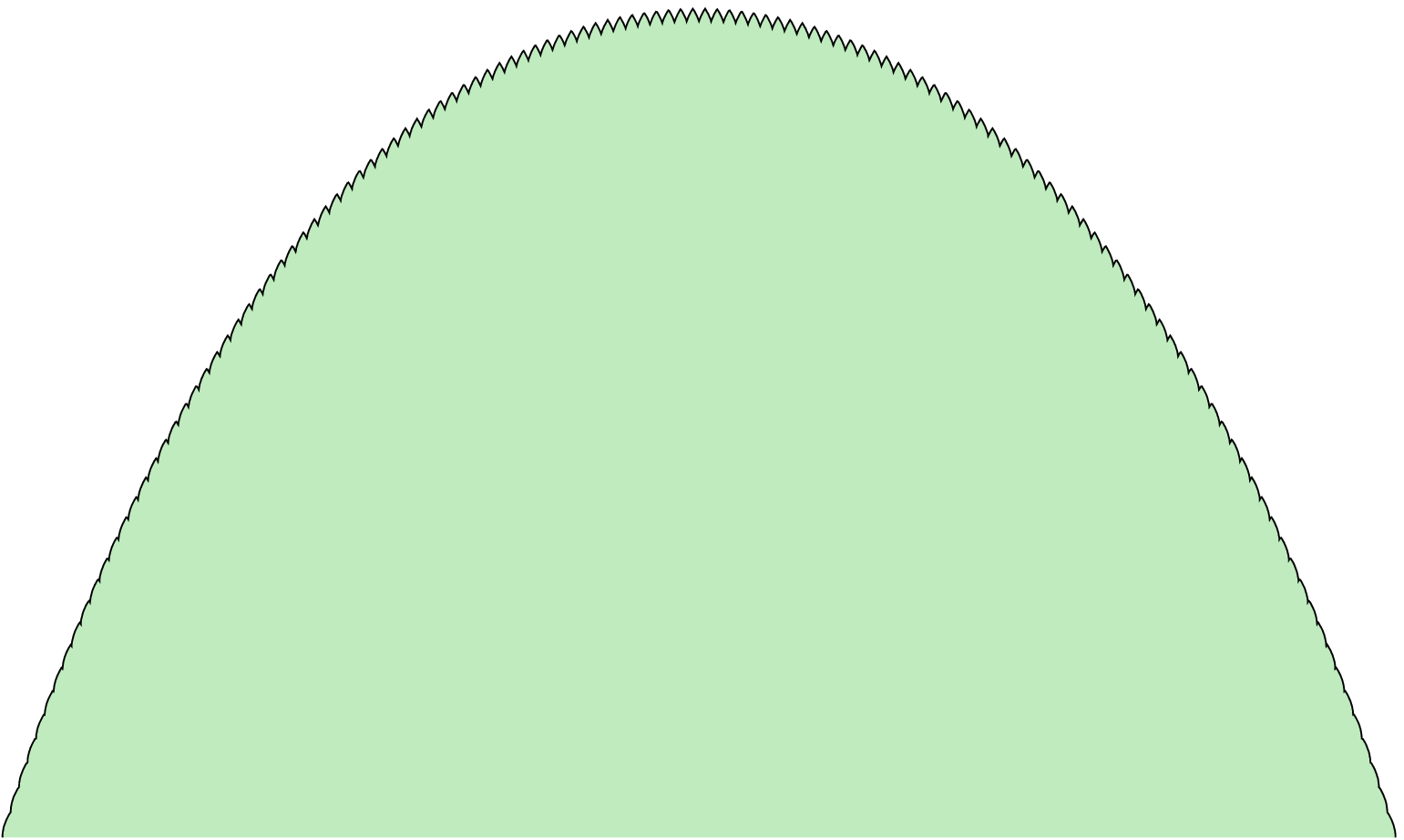}
   \caption{Limiting curves observed for the generalized Pascal-adic transformations:
   $d=3$ (left), $d=8$ (middle) and $d=128$ (right).}
   \label{generalized}
\end{figure}

In \cite{mela04}, Xavier Méla introduced a family of transformations generalizing the 
Pascal-adic transformation. They can be constructed following the same cutting and stacking procedure as described in Appendix~\ref{constructionPA}, but in which each tower is split into $d$ sub-columns; the last $(d-1)$ sub-columns of the tower $\tau_{n,k}$ being sent to the first $(d-1)$ sub-columns of the tower $\tau_{n,k+1}$. 
(The standard Pascal-adic transformation corresponds to the particular case $d=2$.)
Numerical simulations (see Figure~\ref{generalized}) indicate that similar results of convergence as those proved in the present paper also hold in this more general 
context. The limiting curves also seem to be self-affine, but defined with $d$ affinities instead of just 2. Interestingly, as $d\to\infty$ these curves seem to 
converge to a smooth function.

\subsubsection{Rotations and rank-one transformation}

Many questions remain open concerning the Pascal-adic transformation. Its mixing
properties are totally unknown, but it is conjectured that it is at least weakly mixing.
Related to this important question, we can ask whether such behaviour of ergodic sums can
be observed in systems defined by an irrational rotation on the circle.

One of the few properties which have been established for the Pascal-adic transformation
is the loose-Bernoullicity (see~\cite{Janvresse-delaRue04}). In the class of zero-entropy
systems, to which belongs the Pascal-adic, loose-Bernoullicity is the weaker of a
sequence of ergodic properties:
\begin{quotation}
rank one $\Longrightarrow$ finite rank $\Longrightarrow$ local rank one $\Longrightarrow$ loosely-Bernoulli.
\end{quotation} 
Méla and Petersen ask in~\cite{mela-petersen04} whether those stronger properties are satisfied by
the Pascal-adic transformation. The conjecture is that it is not even of local rank one.
Connected to this problem, it would be interesting to study the behaviour of the corrections
to the ergodic theorem in general rank-one systems. Can phenomenon such as those established
in this work appear in the rank-one category?
% More generally, what can be said about such corrections to the ergodic 
% theorem for rank-one systems? This is to be connected with this other question
% on the Pascal-adic transformation: is this transformation of rank one? Just recall
% that the rank-one property implies the loose-Bernoullicity, which has been established
% for this system in~\cite{Janvresse-delaRue04}.

\medskip

\appendix
\section{Construction of the Pascal-adic transformation}
\label{constructionPA}
\begin{figure}[t!] 
  \label{construction}
  \centering
\begin{picture}(0,0)%
\includegraphics{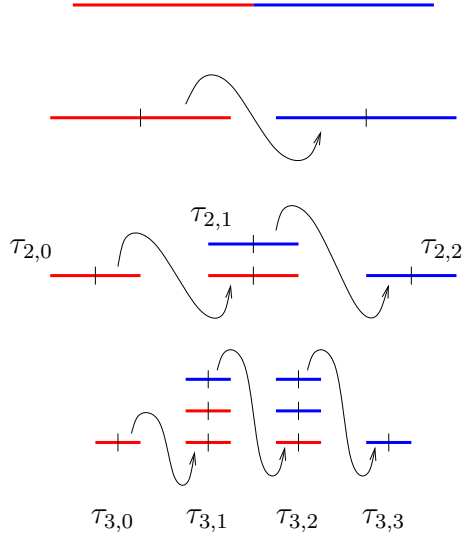}%
\end{picture}%
\setlength{\unitlength}{2486sp}%
\begingroup\makeatletter\ifx\SetFigFont\undefined%
\gdef\SetFigFont#1#2#3#4#5{%
  \reset@font\fontsize{#1}{#2pt}%
  \fontfamily{#3}\fontseries{#4}\fontshape{#5}%
  \selectfont}%
\fi\endgroup%
\begin{picture}(4488,5289)(-629,-6892)
\put(181,-6811){\makebox(0,0)[lb]{\smash{{\SetFigFont{12}{14.4}{\rmdefault}{\mddefault}{\updefault}{\color[rgb]{0,0,0}$\tau_{3,0}$}%
}}}}
\put(2026,-6811){\makebox(0,0)[lb]{\smash{{\SetFigFont{12}{14.4}{\rmdefault}{\mddefault}{\updefault}{\color[rgb]{0,0,0}$\tau_{3,2}$}%
}}}}
\put(1126,-6811){\makebox(0,0)[lb]{\smash{{\SetFigFont{12}{14.4}{\rmdefault}{\mddefault}{\updefault}{\color[rgb]{0,0,0}$\tau_{3,1}$}%
}}}}
\put(2881,-6811){\makebox(0,0)[lb]{\smash{{\SetFigFont{12}{14.4}{\rmdefault}{\mddefault}{\updefault}{\color[rgb]{0,0,0}$\tau_{3,3}$}%
}}}}
\put(-629,-4111){\makebox(0,0)[lb]{\smash{{\SetFigFont{12}{14.4}{\rmdefault}{\mddefault}{\updefault}{\color[rgb]{0,0,0}$\tau_{2,0}$}%
}}}}
\put(3466,-4111){\makebox(0,0)[lb]{\smash{{\SetFigFont{12}{14.4}{\rmdefault}{\mddefault}{\updefault}{\color[rgb]{0,0,0}$\tau_{2,2}$}%
}}}}
\put(1171,-3796){\makebox(0,0)[lb]{\smash{{\SetFigFont{12}{14.4}{\rmdefault}{\mddefault}{\updefault}{\color[rgb]{0,0,0}$\tau_{2,1}$}%
}}}}
\end{picture}%
  \caption{Cutting and stacking construction of the Pascal-adic transformation}
  \end{figure}
Here we recall the construction  of the Pascal-adic transformation,
following the cutting and stacking model exposed in \cite{mela-petersen04}.
Our space $X$ is the interval $[0,1[$, equipped with its Borel $\sigma$-algebra $\mathscr{A}$. 

We start by dividing $X$ into two subintervals
$P_0:= [0,1/2[$ and $P_1:= [1/2,1[$. Let $\P_1:=\{P_0,P_1\}$ be the partition obtained at this first step. 
We also consider $P_0$ and $P_1$ as ``degenerate'' Rokhlin
towers of height 1, respectively denoted by $\tau_{1,0}$ and $\tau_{1,1}$.

On second step, $P_0$ and $P_1$ are divided into two equal subintervals. The transformation $T$ is defined
on the right piece of $P_0$ by sending it linearly onto the left piece of $P_1$. This gives a collection of 3 disjoint Rokhlin towers %(stacks) 
denoted by $\tau_{2,0}, \tau_{2,1}$, $\tau_{2,2}$, with respective heights 1, 2, 1 (see figure~\ref{construction}).

After step $n$, we get $(n+1)$ towers $\tau_{n,0},\ldots,\tau_{n,n}$, with respective heights $\binom{n}{0},\ldots,\binom{n}{n}$,
the width of $\tau_{n,k}$ being $2^{-n}$. 
At this step, the transformation $T$ is defined on the whole space except
the top of each stack. We then divide each stack into two sub-columns with equal width, and define $T$ on the right piece of the top
of $\tau_{n,k}$ by sending it linearly onto the left piece of the base  of $\tau_{n,k+1}$.
Repeating recursively this construction, $T$ is finally defined on all of $X$ 
except on countably many points. 
It is well known that the ergodic invariant measures for this transformation
are given by the one-parameter family $(\mu_p)_{0<p<1}$, where $\mu_p$ is the image of the Bernoulli measure $B(1-p,p)$ on $\{0,1\}^{\mathbb{N}}$ by the application $(x_k)\longmapsto\sum_{k\ge1}x_k/2^k$. This measure
$\mu_p$ can be interpreted as follows: For each $x\in[0,1[$ and $n\ge1$, denote by $k_n(x)$ 
the unique index such that 
$x\in\tau_{n,k_n(x)}$. Under $\mu_p$, conditionned on $k_1(x),k_2(x),\ldots,k_n(x)$, the value of
$k_{n+1}(x)$ is either $k_n(x)$ (with probability $1-p$) or $k_{n}(x)+1$  (with probability $p$). Thus, the
law of large number gives
\begin{equation}
\dfrac{k_n(x)}{n}\ \mathop{\longrightarrow}_{n\rightarrow\infty}\ p\qquad\mu_p-\text{a.e.}
\end{equation}

\section{Links with Conway recursive sequence}
\label{app_conway}
Let us recall that Conway recursive sequence is defined by $C(1)=C(2)=1$, and for $j\ge3$
$$ C(j)\ =\ C(C(j-1))+C(j-C(j-1)). $$
One easily checks that the differences
$$ \Delta C(j)\ :=\ C(j)-C(j-1) $$
are always 0 or 1. Following Mallows \cite{Mallows1991}, it is convenient to introduce the sequence
$$ D(j)\ :=\ 2\Delta C(j)-1\ \in\ \{-1,1\}.  $$
It is shown in \cite{Mallows1992} that the sequence $D(j)$, $j\ge3$ is obtained by the concatenation
of the $B_{n,k}$'s after substituting $a$ by 1 and $b$ by $-1$:
$$ (D(j))_{j\ge3}\ =\ B_{1,0}B_{1,1}B_{2,0}B_{2,1}B_{2,2}B_{3,0}\dots $$
This is a consequence of the following remarkable property: Recall that $B_{n,k}$ is the concatenation
of $B_{n-1,k-1}$ and $B_{n-1,k}$; in fact  $B_{n,k}$ can also be obtained by this alternative procedure.
Cut $B_{n-1,k-1}$ after each $a$, and $B_{n-1,k}$ after each $b$. Interleaving the resulting pieces
produces $B_{n,k}$. For example $B_{4,2}=aababb$ is the concatenation of $a|a|b$ and $ab|b|$ and can be written
as $a|ab|a|b|b$.

The graph of the function with increments $D(j)$ consists in a series of humps corresponding to 
intervals $2^n+1\le j\le 2^{n+1}$. In most works dealing with Conway sequence, the asymptotic shape of these humps
is studied, and it is shown to be given by some smooth explicit function. Each of these humps corresponds
to the graph associated to the concatenation of all the words $B_{n,k}$ on a given line. It turns out that at this
scale, the fractal structure is lost. The results presented in our paper can thus also be interpreted as the analysis of small fluctuations in the convergence of these humps.

\bibliographystyle{amsplain}

\providecommand{\bysame}{\leavevmode\hbox to3em{\hrulefill}\thinspace}
\providecommand{\MR}{\relax\ifhmode\unskip\space\fi MR }
\providecommand{\MRhref}[2]{%
  \href{http://www.ams.org/mathscinet-getitem?mr=#1}{#2}
}
\providecommand{\href}[2]{#2}

%\bibliography{macdo}

\end{document}